\documentclass[10pt]{amsart}
\usepackage{amscd} 
\usepackage{amsfonts} 
\usepackage{amssymb} 
\usepackage{latexsym}

\newcommand{\ncm}{\newcommand}

\newtheorem{theorem}{Theorem}[section]
\newtheorem{prop}[theorem]{Proposition}
\newtheorem{lemma}[theorem]{Lemma}

\newtheorem{lem&def}[theorem]{Lemma \& Definition}
\newtheorem{definition}[theorem]{Definition}
\newtheorem{example}[theorem]{Example}

\def\C{\mathbb{C}\,}

\def\id{\mbox{\rm id}}

\def\into{\hookrightarrow}
\def\to{\rightarrow}

\def\End{\mbox{\rm End}\,}

\def\Hom{\mbox{\rm Hom}\,}
\def\o{{\otimes}}    
\def\|{\, | \, }
\def\bra{\langle}
\def\ket{\rangle}

\ncm{\rarr}[1]{\stackrel{#1}{\longrightarrow}}
\ncm{\larr}[1]{\stackrel{#1}{\longleftarrow}}

\def\cop{\Delta}

\def\eps{\varepsilon}

\def\du1{\hat 1}

\def\-1{_{(-1)}}
\def\0{_{(0)}}
\def\1{_{(1)}}
\def\2{_{(2)}}
\def\3{_{(3)}}

\def\du1{\hat 1}

\def\ract{\triangleleft}

\begin{document}

\title[Depth three towers]{Depth three towers of rings and groups}
\author{Lars Kadison} 
\address{Department of Mathematics \\ University of Pennsylvania \\
David Rittenhouse Lab, 209 S. 33rd St. \\ 
Philadelphia, PA 19104} 
\email{lkadison@math.upenn.edu}
\thanks{My thanks to David Harbater and the Penn Galois seminar
for the invitation.}
\date{} 

\begin{abstract}
Depth three and finite depth are notions known for subfactors via diagrams and Frobenius extensions of rings via centralizers in endomorphism towers.  From the point of view of depth two ring extensions, we provide a clear definition of
depth three for a tower of three rings $C \subseteq B \subseteq A$. If $A = \End B_C$ and $B \| C$
is a Frobenius extension, this captures the notion of depth three for a Frobenius extension.  For example we provide
an algebraic proof that if $B \| C$ is depth three, then $A \| C$
is depth two.  
If $A$, $B$ and $C$ correspond to a tower of subgroups $G > H > K$ via 
the group algebra over a fixed base ring, the depth three condition
is the condition that subgroup $K$ has normal closure $K^G$ contained in $H$.
For a depth three tower of rings, there is  an interesting algebraic theory from the point of view of Galois correspondence
theory for the ring $\End {}_BA_C$ and coring $(A \o_B A)^C$ with respect to the centralizers $A^B$ and $A^C$ involving Morita context bimodules, nondegenerate pairings and right coideal subrings.   
\end{abstract} 
\maketitle

\section{Introduction}

Depth $n$ is a notion from the classification of subfactors which describes where
in the derived tower of centralizers, if at all,  there occurs three successive algebras forming
a basic construction $C \into B \into \End B_C$.
Depth two plays the most important role in finite depth classification theory \cite{NV}. This is partly because
a finite depth subfactor embeds via
its Jones tower into a depth two subfactor
(see Theorem~\ref{prop-d2d3} for the
depth three algebraic version).  
A subfactor $B \subseteq A$ is depth two then if the centralizers $V_A(B) \into V_{A_1}(B) \into
V_{A_2}(B)$ is a basic construction, where $A \into A_1 \into A_2 \into A_3$ is a Jones tower
of iterated basic constructions.  The subfactor $B \subseteq A$ is depth three if
instead the centralizers $V_{A_1}(B) \into V_{A_2}(B) \into V_{A_3}(B)$ is 
a basic construction.
Th algebraic property of finite depth may be descibed most easily starting with a Frobenius extension
$A \supseteq B$, where the definition guarantees the existence of a bimodule homomorphism
$A \to B$ with dual bases for the finitely generated projective $B$-module $A$ \cite{KN}.

A careful algebraic study of the depth two condition on subalgebra
$B \subseteq A$ shows that it is most simply expressed
as a type of central projectivity condition on the tensor-square $A \o_B A$
w.r.t.\  $A$ as natural $A$-$B$-bimodules and $B$-$A$-bimodules.  There is 
a Galois theory connected to this viewpoint with Galois quantum groupoids, 
in the category of Hopf algebroids. Although a future viewpoint on
 depth two ring extension in this generality
might be that it is better called a ``normal  extension,''
depth two does presently suggest that it is part of a larger theory of depth $2, \, 3$ and beyond
 for ring extensions.  
Indeed depth three does lend itself, after reformulation, to
a notion for ring extensions (see
the preprint version of \cite{KS}).
However, in this preprint we prefer to
view depth three 
as a property most naturally associated to a tower of three algebras or rings
$C \subseteq B \subseteq A$.  This tower is right depth three (rD3) if $A \o_B A$
is $A$-$C$-isomorphic to a direct summand of $A \oplus \cdots \oplus A$.  
The advantage of this definition over the one in \cite[preprint version]{KS}
is that it is  close to the depth two definition so that a  substantial
amount of depth two theory is available as we see in this paper.  At the same
time, it seems inevitable that depth three
will play a role in any Galois correspondence theory involving depth
two, simple algebras and purely
algebraic hypotheses: see Theorem~\ref{th-conv}, the last
section and compare with \cite{NV}.  
The connection with classical depth three subfactors may be seen as follows:
if $C \subseteq B$ is a Frobenius extension with $A = \End B_C$, it
follows that $A \cong B \o_C B$, that ${}_AA \o_B A_C$ reduces to ${}_AB \o_C B \o_C B_C$
and ${}_AA_C$ to ${}_AB \o_C B_C$, the terms
in which the depth three condition is expressed  in \cite[preprint version]{KS}.  

In section 2 of this paper we note that right or left D3 ring towers are characterized in terms either
of the tensor-square, H-equivalent modules, quasibases or the endomorphism ring. In section~3
we show that a tower of subgroups $G > H > K$ of finite index with the condition that the
normal closure $K^G < H$ ensures that the group algebras $ F[G] \supseteq F[H] \supseteq F[K]$
are a depth three tower w.r.t.\ any base ring $F$.  We propose that the converse
is true if $G$ is a finite group and $F = \C$.  In section~4 we study the right
coideal subring $E = \End {}_BA_C$
as well as the bimodule and co-ring $P = (A \o_B A)^C$, which provide the quasibases for a right D3 tower $A \| B \| C$.  
We show that right depth three towers may be characterized by $P$ being finite projective as
a left module over the centralizer $V = A^C$ and a pre-Galois isomorphism $A \o_B A \stackrel{\cong}{\longrightarrow} A \o_V P$.

\section{Definition and first properties of depth three towers}

Let $A$, $B$ and $C$ denote rings with identity element,
and $C \to B$, $B \to A$ denote ring homomorphisms preserving the
identities.  We use ring extension notation $A \| B \| C$ for
$ C \to B \to A$ and call this a tower of rings: an important
special case if of course $C \subseteq B \subseteq A$ of subrings
$B$ in $A$ and $C$ in $B$.  Of most importance to us are the induced
bimodules such as ${}_BA_C$ and ${}_CA_B$.  We may naturally also
choose to work with algebras over commutative rings, and obtain almost
identical results.

We denote the centralizer subgroup of a ring $A$ in an $A$-$A$-bimodule $M$ by
$M^A =\{ m \in M \| \, \forall a \in A, ma = am \}$.  We also
use the notation $V_A(C) = A^C$ for the centralizer subring of $C$ in $A$.
This should not be confused with our notation $K^G$ for the normal
closure of a subgroup $K < G$.  Notation like $\End B_C$
will denote the ring of endomorphisms of the module $B_C$
under composition and  addition (and not algebra homomorphisms fixing $C$ or the like).       

\begin{definition}

A tower of rings $A \| B \| C$ is right depth three (rD3) if the tensor-square
$A \o_B A$ is isomorphic as $A$-$C$-bimodules to a direct summand of
 a finite direct
sum of $A$ with itself: in module-theoretic symbols, this becomes, for some positive integer $N$, 
\begin{equation}
\label{eq: rD3}
{}_A A\o_B A_C \oplus * \cong {}_A A^N_C
\end{equation}
\end{definition}
By switching to $C$-$A$-bimodules instead, we naturally define
a left D3 tower of rings. The theory for these is entirely analogous
to rD3 towers and is briefly considered at the end of this section.

Recall that over a ring $R$, two modules $M_R$ and $N_R$ are
H-equivalent if $M_R \oplus * \cong N^n_R$ and $N_R \oplus * \cong M_R^m$
for some positive integers $n$ and $m$.  In this case, the endomorphism
rings $\End M_R$ and $\End N_R$ are Morita equivalent with
context bimodules $\Hom (M_R, N_R)$ and $\Hom (N_R, M_R)$.

\begin{lemma}
A tower $A \| B \| C$ of rings is rD3 iff the natural $A$-$C$-bimodules
$A \o_B A$ and $A$ are H-equivalent.
\end{lemma}
\begin{proof}  
We note that for any tower of rings, $A \oplus * \cong A \o_B A$ as $A$-$C$-bimodules, since the epi $\mu: A \o_B A \to A$ splits as an
$A$-$C$-bimodule arrow.  
\end{proof}

Since for any tower of rings $\End {}_AA_C$ is isomorphic to the centralizer
$V_A(C) = A^C$ (or anti-isomorphic according to convention),
we see from the lemma that the notion of rD3 has something to do with
classical depth three. Indeed,

\begin{example}
\begin{rm}
If $B \| C$ is a Frobenius extension, with Frobenius system $(E,x_i,y_i)$
satisfying for each $a \in A$,  
\begin{equation}
\sum_i E(a x_i)y_i = a = \sum_i x_i E(y_i a)
\end{equation}
then $ B \o_C B \cong \End B_C := A$ via $x \o_B y \mapsto \lambda_x \circ E \circ \lambda_y$ for left multiplication $\lambda_x$ by element $x \in B$.
Let $B \to A$ be this mapping $B \into \End B_C$ given by $b \mapsto \lambda_b$.
It is then easy to show that ${}_AB \o_C B \o_C B_C \cong {}_AA \o_B A_C$,
so that condition~(\ref{eq: rD3}) is equivalent to the condition for rD3 in preprint \cite{KS}, which in turn slightly generalizes the condition in \cite{KN}
for depth three free Frobenius extension.  We should make note here that right or left depth three would be equivalent notions for Frobenius extensions,
since $\End B_C$ and $\End {}_CB$ are anti-isomorphic for such.  
\end{rm}
\end{example}

Another litmus test for a correct notion of depth three is that depth two
extensions should be depth three
in a certain sense. Recall that
a ring extension $A \| B$ is right depth two (rD2)
if the tensor-square $A \o_B A$ is $A$-$B$-bimodule
isomorphic to $N$ copies of $A$ in a direct sum with itself:
\begin{equation}
\label{eq: rD2}
{}_AA \o_B A_B \oplus * \cong {}_AA^N_B
\end{equation}
 Since the notions pass from ring extension
to tower of rings, there are several cases to look at.  

\begin{prop}
\label{prop-d2d3}
Suppose $A \| B \| C$ is a tower of rings. We note:
\begin{enumerate}
\item If $B= C$ and $B \to C$ is the identity mapping,
then $A \| B \| C$ is rD3 $\Leftrightarrow$ $A \| B$ is rD2. 
\item If $A \| B$ is rD2, then $A \| B \| C$ is rD3 w.r.t.\ any ring extension
$B \| C$.
\item If $A \| C$ is rD2 and $B \| C$ is a separable extension,
then $A \| B \| C$ is rD3.    
\item If $B \|C$ is left D2, and $A = \End B_C$, then 
$A \| B \| C$ is left D3.
\item If $C$ is the trivial subring, any ring extension $A \| B$,
where ${}_BA$ is finite projective, together
with $C$ is rD3.  
\end{enumerate}
\end{prop}
\begin{proof}
The proof follows from comparing eqs.~(\ref{eq: rD3}) and~(\ref{eq: rD2}),
noting that $A \o_B A \oplus * \cong A \o_C A$ as natural $A$-$A$-bimodules
if $B \| C$ is a separable extension (and having a separability element 
$e = e^1 \o_C e^2 \in (B \o_C B)^B$ satisfying $e^1 e^2 = 1$),
and finally from \cite{LK2006A} that $B \| C$ left D2 extension $\Rightarrow$  $A \| B$ is left D2 extension if $A = \End B_C$. The last
statement follows from tensoring ${}_BA \oplus * \cong {}_BB^n$ by
${}_AA\o_B -$.  
\end{proof}

The next theorem is a converse and algebraic
simplification of a key fact in subfactor Galois theory (the n = 3 case):  a depth three subfactor $N \subseteq M$ yields a depth
two subfactor $N \subseteq M_1$, w.r.t.\
its basic constuction $M_1 \cong M \o_N M$. In preparation,
let us call a ring extension $B \| C$ rD3 if the endomorphism
ring tower $A \| B \| C$ is rD3, where $A = \End B_C$ and
$A \| B$ has underlying map $\lambda:
B \rightarrow \End B_C$, the left regular mapping
 given by $\lambda(x)(b) = xb$
for all $x, b \in B$.

\begin{theorem}
\label{th-conv}
Suppose $B \| C$ is a Frobenius extension
and $A = \End B_C$.  If $ B \| C$ is
rD3, then the composite extension $A \| C$ is D2. 
\end{theorem}
\begin{proof}
There is a well-known bimodule isomorphism for a Frobenius
extension $B \| C$,  between its endomorphism ring
and  its tensor-square, 
${}_BA_B \cong  {}_B B \otimes_C B_B$.
Tensoring by ${}_AA \o_B - \o_B A_A$,
we obtain $A \o_C A \cong A \o_B A \o_B A$
as natural $A$-$A$-bimodules.  Now
restrict the bimodule isomorphism in eq.~(\ref{eq: rD3}) on the right to $B$-modules and tensor by ${}_AA \o_B -$
to obtain ${}_AA \o_C A_C \oplus * \cong
{}_A A \o_B A_C^N$ after substitution of
the tensor-cube over $B$ by the tensor-square
over $C$.  By another application of  eq.~(\ref{eq: rD3})
we arrive at $${}_AA \o_C A_C \oplus * \cong
{}_A A^{N^2}_C$$
Thus $A \| C$ is right D2.  Since it is
a Frobenius extension as well, it is also
left depth two.  
\end{proof}

We introduce quasi-bases for right depth three towers.

\begin{theorem}
\label{th-qb}
A tower $A \| B \| C$ is right depth three iff there are $N$ elements each
of $\gamma_i \in \End {}_BA_C$ and of $u_i \in (A \o_B A)^C$ satisfying
(for each $x, y \in A$)
\begin{equation}
\label{eq: rd3qb}
x \o_B y = \sum_{i=1}^N x \gamma_i(y) u_i
\end{equation}
\end{theorem}
\begin{proof}
From the condition~(\ref{eq: rD3}), there are obviously $N$ maps each of 
\begin{equation}
f_i \in \Hom ({}_AA_C, {}_AA \o_B A_C), \ \ g_i \in \Hom ({}_AA \o_B A_C, {}_AA_C)
\end{equation}
such that $\sum_{i = 1}^N f_i \circ g_i = \id_{A \o_B A}$.
First, we note that for any tower of rings, not necessarily rD3, 
\begin{equation}
\Hom ({}_AA_C, {}_AA \o_B A_C) \cong (A \o_B A)^C
\end{equation}
via $f \mapsto f(1_A)$.  The inverse is given by $p \mapsto ap$
where $p = p^1 \o_B p^2 \in (A \o_B A)^C$ using a Sweedler-type notation
that suppresses a possible summation over simple tensors.

The other hom-group above also has a simplification.  We note
that for any tower,
\begin{equation}
\Hom ( {}_AA \o_B A_C, {}_AA_C) \cong \End {}_BA_C
\end{equation}
via $F \mapsto F(1_A \o_B -)$.  Given $\alpha \in \End {}_BA_C$,
we define an inverse sending $\alpha$ to the homomorphism
$x \o_B y \mapsto x \alpha(y))$.  

Let $f_i$ correspond to $u_i \in (A \o_B A)^C$ 
and $g_i$ correspond to $\gamma_i \in \End {}_BA_C$ via the mappings
just described.  We compute:
$$ x \o_B y = \sum_i f_i (g_i(x \o y)) = \sum_i f_i (x \gamma_i(y))=
\sum_i x \gamma_i(y)u_i, $$
which establishes the rD3 quasibases equation in the theorem,
given an rD3 tower.

For the converse, suppose we have $u_i \in (A \o_B A)^C$ and
$\gamma_i \in \End {}_BA_C$ satisfying the equation in the theorem.
 Then map $\pi: A^N \to A \o_B A$ by $$\pi: (a_1,\ldots,a_N) \longmapsto
\sum_i a_iu_i,$$ an $A$-$C$-bimodule epimorphism split
by the mapping $\sigma: A \o_B A \into A^N$ given by
$$\sigma(  x \o_B y ) := (x\gamma_1(y),\ldots,x\gamma_N(y)). $$
It follows from the equation above
 that $\pi \circ \sigma = \id_{A \o_B A}$.
\end{proof}

\subsection{Left D3 towers and quasibases} A tower of rings $A \| B \| C$
is left D3 if the tensor-square $A \o_B A$ is an $C$-$A$-bimodule direct
summand of $A^N$ for some $N$.  If $B = C$, this recovers the definition
of a left depth two extension $A \| B$.  There is a left version of
all results in this paper:  we note that  $A \| B \| C$ is a right D3
tower if and only if $A^{\rm op} \| B^{\rm op} \| C^{\rm op}$ is
a left D3 tower (cf.\ \cite{LK2006A}). 

The next theorem refers to notation established in the example above.  
\begin{theorem}
Suppose $B \| C$ is a Frobenius extension with $A = \End B_C$.  
 Then $A \| B \| C$ is right depth three if and only if 
$A \| B \| C$ is left depth three.
\end{theorem}
\begin{proof}
It is well-known that also $A \| B$ is a Frobenius extension. 
Then $A \o_B A \cong \End A_B$ as natural $A$-$A$-bimodules.
Also $A \o_B A \cong \End {}_BA$ by a similar mapping utilizing
the Frobenius homomorphism in one direction, and the dual bases
in the other.  As a result, the left and right endomorphism rings
are anti-isomorphic.  

Now note the following characterization of left D3 with proof
almost identical with that of \cite[Prop.\ 3.8]{LK2007}:
If $A \| B \| C$ is a tower where $A_B$ if finite projective, then
$A \| B \| C$ is left D3 $\Leftrightarrow$ $\End A_B \oplus * \cong A^N$
as natural $A$-$C$-bimodules.
The proof involves noting that $\End A_B \cong \Hom (A \o_B A_A,A_A)$
via $$f \longmapsto (a \o a' \mapsto f(a)a').$$

Similarly, if $A \| B \| C$ is a tower where ${}_BA$ is finite projective,
then $A \| B \| C$ is right D3 if and only if $\End {}_BA \oplus * \cong A^N$
as natural $C$-$A$-bimodules.    
 
Of course a Frobenius extension satisfies both finite projectivity conditions.
The anti-isomorphism of the left and right endomorphism rings twists
the $C$-$A$-structure to an $A$-$C$-structure, thereby demonstrating
the equivalence of left and right D3 conditions on $A \o_B A$ relative to  $A \cong \End A_A$.  
\end{proof}

In a fairly obvious reversal to
 opposite ring structures in the proof of theorem~\ref{th-qb},
we see that a tower $A \| B \| C$ is left D3 iff there are
$N$ elements $\beta_j \in \End {}_CA_B$ and $N$ elements $t_j \in (A \o_B A)^C$
such that for all $x, y \in A$, we have
\begin{equation}
\label{eq: ld3qb}
x \o_B y = \sum_{j=1}^N t_j \beta_j(x)y
\end{equation}
 
We record the characterization of left D3, noted above in the proof, for towers satisfying a finite
projectivity condition.

\begin{prop}
Suppose $A \| B \| C$ is a tower of rings where $A_B$ is finite projective.
Then this tower is left D3 if and only if the natural $A$-$C$-bimodules satisfy for some $N$, 
\begin{equation}
\End A_B \, \oplus \,  * \, \cong A^N
\end{equation}
\end{prop}

Finally we define a tower $A \| B \| C$ to be D3 if it is both left D3 and right D3.  


\section{Depth three for towers of groups}

Fix a base ring $F$.  Groups give rise to rings via $G \mapsto
F[G]$, the functor associating the group algebra $F[G]$ to a group
$G$.  Therefore we can pull back the notion of depth $2$ or $3$
for ring extensions or towers to the category of groups (so long as
reference is made to the base ring). 

In the paper \cite{KK}, a depth two subgroup w.r.t.\ the complex numbers
is shown to be equivalent to the notion of normal subgroup for finite
groups.  This consists of two results.  The easier result is that over any
base ring, a normal subgroup of finite index is depth two
by exhibiting left or right D2 quasibases via  coset representatives
and projection onto cosets.  This proof suggests that the converse hold
as well.  The second result is a converse for complex finite dimensional D2
group algebras where normality of the subgroup is established using character theory and Mackey's subgroup theorem.

In this section, we will similarly do the first step in showing
what group-theoretic notion corresponds to depth three tower of rings.
Let $G > H > K$ be a tower of groups, where $G$ is a finite group,
$H$ is a subgroup, and $K$ is a subgroup of $H$.  Let $A = F[G]$,
$B = F[H]$ and $C = F[K]$.  Then $A \| B \| C$ is a tower of rings,
and we may ask what group-theoretic notion on $G > H > K$ will guarantee,
with fewest possible hypotheses, that $A \| B \| C$ is rD3.

\begin{theorem}
The tower of groups algebras $A \| B \| C$ is D3 if
the corresponding tower of groups $G > H > K$ satisfies
\begin{equation}
K^G < H
\end{equation}
where $K^G$ denotes the normal closure of $K$ in $G$.  
\end{theorem}
\begin{proof}
Let $\{ g_1,\ldots,g_N \}$ be double coset representatives such
that $G = \coprod_{i=1}^N Hg_iK$.  Define $\gamma_i(g) = 0$
if $g \not \in Hg_i K$ and $\gamma_i (g) = g$ if $g \in Hg_iK$.  
Of course, $\gamma_i \in \End {}_BA_C$ for $i = 1,\ldots, N$.  

Since $K^G \subseteq H$, we have $gK \subseteq Hg$ for each $g \in G$.
  Hence for each $k \in K$, $g_jk = h g_j$ for some $h \in H$.  It follows
that $$g_j^{-1} \o_B g_jk = g_j^{-1} h \o_B g_j = kg_j^{-1} \o_B g_j.$$

Given $g \in G$, we have $g = hg_j k$ for some $j = 1,\ldots,N$,
$h \in H$, and $k \in K$.
Then we compute:
$$ 1 \o_B g = 1 \o_B hg_jk = hg_j g_j^{-1} \o_B g_jk = hg_jk g_j^{-1} \o_B g_j $$
so $1 \o_B g =  \sum_i \gamma_i(g) g_i^{-1} \o_B g_i $
where $g_i^{-1} \o_B g_i \in (A \o_B A)^C$.  By theorem then,
$A \| B \| C$ is an rD3 tower.

The proof that the tower of group algebras is left D3 is entirely symmetical
via the inverse mapping.    
\end{proof}

The theorem is also valid for infinite groups where the index $[G : H]$
is finite, since $HgK = Hg$ for each $g \in G$.  

Notice how the equivalent notions of depth two and normality
for finite groups over $\C$ yields the proposition~\ref{prop-d2d3}
for groups.  Suppose we have a tower of groups $G > H > K$ where
$K^G \subseteq H$.  If $K = H$, then $H$ is normal (D2) in $G$.
If $K = \{ e \}$, then it is rD3 together with any subgroup $H < G$.
If $H \ract G$ is a normal subgroup, then necessarily $K^G \subseteq H$.  
If $K \ract G$, then $K^G = K < H$ and the tower is D3.

Question:  Can the character-theoretic proof in \cite{KK} be adapted
to prove that a D3 tower $\C[G] \supseteq \C[H] \supseteq \C[K]$
where $G$ is a finite group 
satisfies $K^G < H$?


\section{Algebraic structure on $\End {}_BA_C$ and $(A \o_B A)^C$}

In this section, we study the calculus of some structures definable for an
rD3 tower $A \| B \| C$, which reduce to the dual bialgebroids over
the centralizer of a ring extension in case $B = C$ and their actions/coactions.
Throughout the section, $A \| B \| C$ will denote a right depth three
tower of rings, 
$$ P := (A \o_B A)^C, \ \ Q := ( A \o_C A)^B,  $$
which are bimodules with respect to the two rings familiar from depth
two theory,
$$ T := (A \o_B A)^B, \ \ U := (A \o_C A)^C$$
Note that $P$ and $Q$ are isomorphic to two $A$-$A$-bimodule Hom-groups: 
\begin{equation}
\label{eqs: tys}
 P \cong \Hom (A \o_C A, A \o_B A), \ \ Q \cong \Hom (A \o_B A, A \o_C A).
\end{equation}
Recall that $T$ and $U$ have multiplications given by
$$ tt' = {t'}^1 t^1 \o_B t^2 {t'}^2, \ \ uu' = {u'}^1u^1 \o_C u^2 {u'}^2, $$
where $1_T = 1_A \o 1_A$ and a similar expression for $1_U$.  
Namely, the bimodule ${}_TP_U$ is given by
\begin{equation}
 {}_TP_U:\ t \cdot p \cdot u = u^1 p^1 t^1 \o_B t^2 p^2 u^2 
\end{equation}
The bimodule ${}_UQ_T$ is given by
\begin{equation}
{}_UQ_T: \ u \cdot q \cdot t = t^1 q^1 u^1 \o_C u^2 q^2 t^2
\end{equation}

We have the following result, also mentioned in passing
in \cite{LK2006B} with several additional hypotheses.

\begin{prop}
\label{prop-ME}
The bimodules $P$ and $Q$ over the rings $T$ and $U$ form
a Morita context with associative multiplications
\begin{equation}
P \o_U Q \rightarrow T,\ \ p \o q \mapsto pq = q^1 p^1 \o_B p^2 q^2 
\end{equation}
\begin{equation}
Q \o_T P \rightarrow U, \ \ q \o p \mapsto qp = p^1 q^1 \o_C q^2 p^2
\end{equation}
If $B \| C$ is an
H-separable extension, then $T$ and $U$ are Morita equivalent rings via this context.
\end{prop}
\begin{proof}
The equations $p(qp') = (pq)p'$ and $q(pq') = (qp)q'$ for
$p,p' \in P$ and $q,q' \in Q$ follow from the four equations directly above.  

Note that
$$ T \cong \End {}_AA \o_B A_A, \ \ U \cong \End {}_AA \o_C A_A $$
as rings.  We now claim that the hypotheses on $A \| B$, $A \| C$ and $B \| C$
imply that the $A$-$A$-bimodules $A \o_B A$ and
$A \o_C A$ are H-equivalent.  Then the endomorphism rings above are Morita
equivalent via context bimodules given by eqs.~(\ref{eqs: tys}), which
proves the proposition.

Since $B \| C$ is H-separable, it is in particular separable, and
the canonical $A$-$A$-epi $A \o_C A \rightarrow A \o_B A$ splits 
via an application of a separability element.  Thus, $A \o_B A \oplus * \cong
A \o_C A$.  
Also, $B \o_C B \oplus * \cong B^N$ as $B$-$B$-bimodules for some positive integer $N$. Therefore, $A \o_C A \oplus * \cong A \o_B A^N$ as $A$-$A$-bimodules
by an application of the functor $A \o_B \, ? \, \o_B A$.  Hence,
$A \o_B A$ and $A \o_C A$ are H-equivalent $A^e$-modules (i.e., $A$-$A$-bimodules). 
\end{proof}
   
We denote the centralizer subrings $A^C$ and $A^B$ of $A$ by
\begin{equation}
 R := V_A(B) \subseteq V_A(C) := V
\end{equation}

We have generalized anchor mappings \cite{LK2006B},
\begin{equation}
R \o_T P \longrightarrow V, \ \ r \o p \longmapsto p^1rp^2
\end{equation}
\begin{equation}
V \o_U Q \longrightarrow R, \ \ v \o q \longmapsto q^1 v q^2
\end{equation}
 
\begin{prop}
The two generalized anchor mappings are bijective if
$B \| C$ is H-separable.
\end{prop}
\begin{proof}
Denote $r \cdot p := p^1 r p^2$ and $v \cdot q := q^1 v q^2$.
From the previous propostion, there are elements $p_i \in P$
and $q_i \in Q$ such that $\sum_i p_i q_i = 1_T$;
in addition, ${p'}_j \in P$ and ${q'}_j \in Q$ such that
$1_U = \sum_j {q'}_j {p'}_j $.  
Let $ v \in V$, then
$$ v = v \cdot 1_U = \sum_j v \cdot ({q'}_j {p'}_j) = \sum_j (v \cdot {q'}_j)\cdot {p'}_j $$
and a similar computation starting with $r = r \cdot 1_T$ shows that the two generalized
anchor mappings are surjective.  

In general, we have the corestriction of the inclusion $T \subseteq A \o_BA$,
 \begin{equation}
{}_TT \into {}_TP
\end{equation}
which is split as a left $T$-module monic by $p \mapsto e^1pe^2$ in case there is a separability element
$e = e^1 \o_C e^2 \in B \o_C B$.  Similarly, 
\begin{equation}
{}_UQ \into {}_UU
\end{equation}
is a split monic in case $B \| C$ is separable.  
Of course, if $B \| C$ is H-separable, we note from  proposition~\ref{prop-ME} and Morita theory 
that $P$ and $Q$ are projective generators on both sides.  

It follows from faithful flatness that the anchor mappings are also injective. 
\end{proof}

Note that $P$ is a $V$-$V$-bimodules (via the commuting homomorphism
and antihomomorphism $V \rightarrow U \leftarrow V$):
\begin{equation}
{}_VP_V:\ \ v \cdot p \cdot v' = vp^1 \o_B p^2 v'
\end{equation}
Note too that $E = \End {}_BA_C$ is an $R$-$V$-bimodule via
\begin{equation}
{}_RE_V: \ \ r \cdot \alpha \cdot v = r\alpha(-)v
\end{equation}
Note the subring and over-ring
\begin{equation}
\End {}_BA_B \subseteq E \subseteq \End {}_CA_C
\end{equation}
which are the total algebras of the left $R$- and $V$- bialgebroids in depth two theory \cite{KS, LK2006A, LK2006B}.  
\begin{lemma}
The modules ${}_VP$ and $E_V$ are finitely generated projective. In case $A \| C$ is
left D2, the subring
$E$ is a right coideal subring of the left $V$-bialgebroid
$\End {}_CA_C$.   
\end{lemma}
\begin{proof}
This follows from eq.~(\ref{eq: rd3qb}), since $p \in P \subseteq A \o_B A$, so 
$$ p = \sum_i p^1 \gamma_i(p^2)u_i $$
where $u_i \in P$ and $p \mapsto p^1 \gamma_i(p^2)$ is in $\Hom ({}_VP, {}_VV)$,
thus dual bases for a finite projective module.  The second claim follows similarly
from
$$ \alpha = \sum_i \gamma_i(-)u_i^1 \alpha(u^2) $$
where $\gamma_i \in E$ and $\alpha \mapsto u^1 \alpha(u^2)$
are mappings in $\Hom (E_V, V_V)$.

Now suppose $\beta_j \in S := \End {}_CA_C$
and $t_j \in (A \o_C A)^C$ are left
D2 quasibases of $A \| C$.  Recall
that the coproduct $\cop: S \rightarrow
S \o_V S$ given by ($\beta \in S$)
\begin{equation}
\cop(\beta) = \sum_j \beta(-t^1_j)t^2_j \o_V \beta_j
\end{equation}
makes $S$ a left $V$-bialgebroid \cite{KS}.
Of course this restricts and corestricts to $\alpha \in E$
as follows:  $\cop(\alpha) \in E \o_V S$.
Hence, $E$ is a right coideal subring of $S$.     
\end{proof}

In fact, if $A \| B$ is also D2,
and $\mathcal{S} = \End {}_BA_B$, 
then $E$ is similarly shown to be 
an $\mathcal{S}$-$S$-bicomodule ring
For we recall the coaction $E \rightarrow
\mathcal{S} \o_R E$ given by 
\begin{equation}
\alpha\-1 \o_R \alpha\0 = \sum_i \tilde{\gamma}_i \o \tilde{u}_i^1 \alpha(\tilde{u}_i^2 -) 
\end{equation}
where $\tilde{\gamma}_i \in \mathcal{S}$
and $\tilde{u}_i \in (A \o_B A)^B$
are right D2 quasibases of $A \| B$
(restriction of \cite[eq.~(19)]{LK2006A}).  
 
Twice above we made use of a $V$-bilinear pairing $P \o E \to V$ given
by
\begin{equation}
\bra p, \alpha \ket := p^1 \alpha(p^2), \ \ (p \in P= (A \o_B A)^C,\, \alpha \in E= \End {}_BA_C)
\end{equation}

\begin{lemma}
The pairing above is nondegenerate.
It induces $E_V \cong \Hom ({}_VP, {}_VV)$ via $\alpha \mapsto \bra - , \alpha \ket$.
\end{lemma}
\begin{proof}
The mapping has the inverse $F \mapsto \sum_i \gamma_i(-)F(u_i)$ where
$\gamma_i \in E, u_i \in P$ are rD3 quasibases for $A \| B \| C$.  
Indeed, $\sum_i \bra p, \gamma_i \ket F(u_i) =$ $ F(\sum_i p^1 \gamma_i(p^2)u_i) = F(p)$
for each $p \in P$ since $F$ is left $V$-linear,
and for each $\alpha \in E$, we note that $\sum_i \gamma_i(-)\bra u_i, \alpha\ket = \alpha$.
\end{proof}

\begin{prop}
There is a $V$-coring structure on $P$ left dual to the ring structure on $E$.
\end{prop}
\begin{proof}
We note that 
\begin{equation}
P \o_V P \cong (A \o_B A \o_B A)^C 
\end{equation}
via $p \o p' \mapsto p^1 \o p^2 {p'}^1 \o {p'}^2$
with inverse $$p = p^1 \o p^2 \o p^3 \mapsto \sum_i (p^1 \o_B p^2 \gamma_i(p^3)) \o_V u_i.$$

Via this identification, define a $V$-linear coproduct $\cop: P \rightarrow P \o_V P$
by \begin{equation}
\cop(p) = p^1 \o_B 1_A \o_B p^2.
\end{equation}
Alternatively, using Sweedler notation and rD3 quasibases,
\begin{equation}
\label{eq: cop}
p\1 \o_V p\2 = \sum_i (p^1 \o_B \gamma_i(p^2)) \o_V u_i 
\end{equation}
   Define a $V$-linear counit $\eps: P \to V$
by $\eps(p) = p^1p^2$.  The counital equations follow readily \cite{BW}.  

Recall from Sweedler \cite{Sw} that the $V$-coring $(P,V,\cop, \eps)$ has
left dual ring ${}^*P := \Hom ({}_VP, {}_VV)$ given by Sweedler notation 
by 
\begin{equation}
(f * g)(p) = f(p\1 g(p\2))
\end{equation}
with $1 = \eps$. Let $\alpha, \beta \in E$.
 If $f = \bra -, \alpha \ket$
and $g = \bra - , \beta \ket$ , we compute $f * g = \bra - , \alpha \circ \beta \ket$
below, which verifies the claim:  
$$ f(p\1 g(p\2)) = \sum_i \bra p^1 \o_B \gamma_i(p^2) \bra u_i, \beta \ket,\, \alpha \ket = 
\bra p^1 \o_B \beta(p^2), \alpha \ket = \bra p, \alpha \circ \beta \ket . $$
\end{proof}

In addition, we note that $P$ is $V$-coring with grouplike element
\begin{equation}
g_P := 1_A \o_B 1_A
\end{equation}
since $\cop(g_P) = 1 \o 1 \o 1 = g_P \o_V g_P$ and $\eps(g_P) = 1$.
  
There is a pre-Galois structure on $A$ given by the right $P$-comodule
structure $\delta: A \rightarrow A \o_V P$, $\delta(a) = a\0 \o_V a\1$
defined by 
\begin{equation}
\delta(a):= \sum_i \gamma_i(a) \o_V u_i.
\end{equation}
The pre-Galois isomorphism $\beta: A \o_BA \stackrel{\cong}{\longrightarrow}
 A \o_V P$ given by \begin{equation}
\beta(a \o a') = a {a'}\0 \o_V {a'}\1
\end{equation}
 is utilized below 
in another characterization of right depth three towers.

\begin{theorem}
A tower of rings $A \| B \| C$ is right depth three if and only if
${}_VP$ is finite projective and $A \o_V P \cong A\o_B A$ as natural
$A$-$C$-bimodules. 
\end{theorem}
\begin{proof}
($\Rightarrow$) If ${}_VP \oplus * \cong {}_VV^N$ and $A \o_V P \cong A \o_B A$,
then tensoring by $A \o_V -$, we obtain $A \o_B A \oplus * \cong A^N$ as natural
$A$-$C$-bimodules, the rD3 defining condition on a tower.

($\Leftarrow$) In proposition~\label{prop-proj} we see that ${}_VP$ is 
f.g.\ projective. Map $A \o_V P \rightarrow A \o_B A$ by $a \o p \mapsto
ap^1 \o_B p^2$, clearly an $A$-$C$-bimodule homomorphism.  
The inverse is the ``pre-Galois'' isomorphism, 
\begin{equation}
\beta: \, A \o_B A \rightarrow A \o_V P, \ \ \beta(a \o_B a') = \sum_i a \gamma_i(a') \o_V  u_i
\end{equation}
 since $\sum_i a p^1 \gamma_i(p^2) \o_V u_i = a \o_V p$ and
$\sum_i a \gamma_i(a') u_i = a \o a'$ for $a,a' \in A, p \in P$.  
\end{proof}

It is an intriguing to continue this study, also for depth $n$ towers,
and study the possibility of an algebraic version of the Galois theory
for subfactors in Nikshych and Vainerman \cite{NV}.  

\end{document}